\documentclass[12pt]{article}
\usepackage{amsmath,amssymb,amsbsy,amsfonts,amsthm,latexsym,
                        amsopn,amstext,amsxtra,euscript,amscd}

\begin{document}

\newtheorem{theorem}{Theorem}
\newtheorem{lemma}[theorem]{Lemma}
\newtheorem{claim}[theorem]{Claim}
\newtheorem{cor}[theorem]{Corollary}
\newtheorem{prop}[theorem]{Proposition}
\newtheorem{definition}{Definition}
\newtheorem{question}[theorem]{Open Question}

\baselineskip 15pt

\title{New estimates of double trigonometric sums with exponential functions}

\author{{M.~Z.~Garaev}
\\
\normalsize{Instituto de Matem{\'a}ticas,  Universidad Nacional
Aut\'onoma de M{\'e}xico}
\\
\normalsize{Campus Morelia, Apartado Postal 61-3 (Xangari)}
\\
\normalsize{C.P. 58089, Morelia, Michoac{\'a}n, M{\'e}xico} \\
\normalsize{\tt garaev@matmor.unam.mx}
\and\\
{A.~A.~Karatsuba}
\\
\normalsize{Steklov Institute of Mathematics}
\\
\normalsize{Russian Academy of Sciences}
\\
\normalsize{GSP-1, ul. Gubkina 8}
\\
\normalsize{Moscow, Russia} \\
\normalsize{\tt karatsuba@mi.ras.ru}
 }


\date{\empty}

\pagenumbering{arabic}

\maketitle


\begin{abstract}
We establish a new bound for the exponential sum
\begin{eqnarray*}
\sum_{x\in\mathcal{X}}\Big|\sum_{y\in \mathcal{Y}}\gamma(y)\exp(2\pi
i a \lambda^{xy}/p)\Big|,
\end{eqnarray*}
where $\lambda$ is an element of the residue ring modulo a large
prime number $p,$ $\mathcal{X}$ and $\mathcal{Y}$ are arbitrary
subsets of the residue ring modulo $p-1$ and $\gamma(n)$ are any
complex numbers with $|\gamma(n)|~\le~1.$ In particular, we improve
several previously known bounds.
\end{abstract}


\section{Introduction}

Let $p$ be a large prime number, $T$ be a divisor of $p-1.$  For a
positive integer $m$ we denote by $\mathbb{Z}_m=\{0, 1,..., m-1\}$
the residue ring modulo $m.$  Let $\lambda$ be an element of
$\mathbb{Z}_p$ of multiplicative order $T,$ i.e.
$\lambda=g^{(p-1)/T}$ for some primitive root $g$ modulo $p.$ Let
$\gamma(n)$ be any complex coefficients with $|\gamma(n)|\le 1.$
Denote
$$
\mathbf{e}_{m}(z)= \exp(2\pi i z/m).
$$

In this paper we investigate the sum
$$
W_{a}(\gamma; T; \mathcal{X}, \mathcal{Y})= \sum_{x\in
\mathcal{X}}\Big|\sum_{y\in\mathcal{Y}}\gamma(y)\mathbf{e}_{p}(a
\lambda^{xy})\Big|,
$$
where $a$ is an integer coprime to $p,$  $\mathcal{X}$ and
$\mathcal{Y}$ are arbitrary subsets of the residue ring
$\mathbb{Z}_{p-1}$ with $|\mathcal{X}|$ and $|\mathcal{Y}|$ elements
correspondingly. This sum is well known and proved to be very
important in many applications, see the recent works~\cite{Ba1}
--\cite{GaKa1}, \cite{Sh} and therein references.

Recently in~\cite{Gar1} a new approach was suggested to estimate
$W_{a}(\gamma; T; \mathcal{X}, \mathcal{Y}).$ One of the advantages
of the new approach compared to the previously known ones is that it
finds its applications also to bound similar sums taken over points
on an elliptic curve, or to bound the corresponding multiplicative
character sums, see~\cite{BFGSh} for the details. The aim of the
present paper is to introduce a new ingredient to the argument
of~\cite{Gar1}  which leads to a new estimate for $W_{a}(\gamma; T;
\mathcal{X}, \mathcal{Y})$ and which, in particular cases, improves
several previously known bounds. Moreover, the new argument in
prospective can also find its applications to modify the
corresponding estimates of~\cite{BFGSh}.

The following statement is the result of our paper.

\begin{theorem}
\label{thm:doublesumcomposit} For any given positive integer $k,$
the following estimate holds:
$$
W_{a}(\gamma; T; \mathcal{X}, \mathcal{Y})\ll
\frac{|\mathcal{X}|^{1-\frac{1}{2k}}|\mathcal{Y}|^{1-\frac{1}{2k+2}}
p^{\frac{1}{2k}+\frac{3}{4k+4}+o(1)}}{T^{\frac{1}{2k+2}}}.
$$
\end{theorem}

Throughout the paper the implied constants in the Landau `$O$' and
`$o$' symbols as well as in the Vinogradov symbols `$\ll$' and
`$\gg$' may depend on the small positive quantity $\varepsilon$ and
on the fixed positive integer $k$. We use $g$ to denote a primitive
root modulo $p.$

\section{Corollaries}

J. B. Friedlander and I. E. Shparlinski~\cite{Fried1} for the sum
$$
W_a(\mathcal{X}, \mathcal{Y})=\sum_{x\in \mathcal{X}}\Big|\sum_{y\in
\mathcal{Y}}\mathbf{e}_{p}(a g^{xy})\Big|
$$
proved the estimate
$$
W_{a}(\mathcal{X}, \mathcal{Y}) \ll
|\mathcal{X}|^{1/2}|\mathcal{Y}|^{5/6}p^{5/8+o(1)}.
$$
If $|\mathcal{X}|\ge |\mathcal{Y}|\ge p^{15/16+\varepsilon},$ then
this estimate provides a nontrivial upper bound for
$S_a(\mathcal{X}, \mathcal{Y}).$ In~\cite{Gar1} the above bound has
been improved to
$$
W_{a}(\gamma; p-1; \mathcal{X}, \mathcal{Y})\ll
|\mathcal{X}|^{1/2}|\mathcal{Y}|^{1/2}p^{7/8+o(1)}
$$
which is nontrivial when $|\mathcal{X}|\ge |\mathcal{Y}|\ge
p^{7/8+\varepsilon}.$ Theorem~\ref{thm:doublesumcomposit} improves
this bound further to the following statement.
\begin{cor}\label{cor:completeset} The estimate
$$
W_{a}(\gamma; p-1; \mathcal{X}, \mathcal{Y})\ll
|\mathcal{X}|^{1-\frac{1}{2k}}|\mathcal{Y}|^{1-\frac{1}{2k+2}}
p^{\frac{1}{2k}+\frac{1}{4k+4}+o(1)}
$$
holds.
\end{cor}
In particular, taking $k$ to be sufficiently large, we see that our
estimate is nontrivial when $|\mathcal{X}|\ge |\mathcal{Y}|\ge
p^{3/4+\varepsilon}.$

Let $\mathcal{X_T}, \mathcal{Y_T}\subset \mathbb{Z}_{T}$ and let
$$
W_a(T; \mathcal{X_T},\mathcal{Y_T})=\sum_{x\in
\mathcal{X_T}}\Big|\sum_{y\in \mathcal{Y_T}}\mathbf{e}_{p}(a
\lambda^{xy})\Big|.
$$
In~\cite{Fried1} it is shown that
$$
|W_a(T; \mathcal{X_T},\mathcal{Y_T})|\ll T^{11/6}p^{1/8}.
$$
In~\cite{Gar1} this bound has been improved to
\begin{equation}
\label{eqn:oldcor} W_a(\gamma; T; \mathcal{X_T},\mathcal{Y_T})\ll
|\mathcal{X_T}|^{1/2}|\mathcal{Y_T}|^{1/2}T^{3/4}p^{1/8+o(1)}.
\end{equation}
From Theorem~\ref{thm:doublesumcomposit} we have the following
consequence.
\begin{cor}
\label{cor:double} For any sets $\mathcal{X_T}\subset\mathbb{Z}_T$
and $\mathcal{Y_T}\subset\mathbb{Z}_T,$ the estimate
\begin{eqnarray*}
W_a(\gamma; T; \mathcal{X_T},\mathcal{Y_T})\ll
|\mathcal{X_T}|^{1-\frac{1}{2k}}|\mathcal{Y_T}|^{1-\frac{1}{2k+2}}T^{\frac{1}{2k}}p^{\frac{1}{4k+4}+o(1)}.
\end{eqnarray*}
holds.
\end{cor}
In particular, if $\mathcal{Y_T}$ is not sufficiently dense in
$\mathbb{Z}_T,$ then Corollary~\ref{cor:double} with $k=1$ improves
the estimate~\eqref{eqn:oldcor}.

To prove Corollary~\ref{cor:double}, we note that in
Theorem~\ref{thm:doublesumcomposit} the sets $\mathcal{X}$ and
$\mathcal{Y}$ are arbitrary subsets of $\mathbb{Z}_{p-1}.$ We shift
the set $\mathcal{X_T}$ by $rT$ with $1\le r\le \frac{p-1}{T}$ and
take $\mathcal{X}$ to be the union of these sets. Analogously we
define $\mathcal{Y}.$ Then
$$
|\mathcal{X}|=\frac{(p-1)|\mathcal{X_T}|}{T},\quad
|\mathcal{Y}|=\frac{(p-1)|\mathcal{Y_T}|}{T}.
$$
Therefore, the estimate of Theorem~\ref{thm:doublesumcomposit} takes
the form
\begin{eqnarray*}
\frac{p^2}{T^2}W_a(\gamma; T; \mathcal{X_T},\mathcal{Y_T})\ll
\left(\frac{p}{T}\right)^{2-\frac{1}{2k}-\frac{1}{2k+2}}\frac{|\mathcal{X_T}|^{1-\frac{1}{2k}}|
\mathcal{Y_T}|^{1-\frac{1}{2k+2}}
p^{\frac{1}{2k}+\frac{3}{4k+4}+o(1)}}{T^{\frac{1}{2k+2}}},
\end{eqnarray*}
whence Corollary~\ref{cor:double}.

It is to be remarked that from Theorem~\ref{thm:doublesumcomposit}
as a consequence one derives the bound
\begin{equation*}
W_{a}(\gamma; T; \mathcal{X},
\mathcal{Y})\ll\frac{|\mathcal{X}|^{1/2}|\mathcal{Y}|^{3/4}p^{7/8+o(1)}}{T^{1/4}}.
\end{equation*}
This bound has been proved in~\cite{GaKa1} in the case when
$\mathcal{Y}$ is an interval. Also note that this bound includes the
best previously known one (see~\cite{Gar1}) and improves it for any
thin set $\mathcal{Y},$ i.e., when $|\mathcal{Y}|\le p^{1-c}$ for
some $c>0.$

\section{The main statement}

For a given divisor $d$ of $p-1,$ let $\mathcal{L}_d$ be any subset
of $\mathbb{Z}_{(p-1)/d}$ such that the elements of $\mathcal{L}_d$
are relatively prime to $(p-1)/d.$ The following Lemma is crucial in
proving Theorem~\ref{thm:doublesumcomposit}.

\begin{lemma}
\label{Lemma} If $d<Tp^{-\frac{1}{2}}(\log p)^{-4k}$ then the
inequality
$$
\sum_{x\in \mathcal{X}}\Big|\sum_{y\in
\mathcal{L}_d}\gamma(dy)\mathbf{e}_{p}(ag^{tdxy})\Big| \ll
\frac{|\mathcal{X}|^{1-\frac{1}{2k}}|\mathcal{L}_d|^{1-\frac{1}{2k+2}}p^{\frac{1}{2k}+\frac{3}{4k+4}}\log
p} {T^{\frac{1}{2k+2}}}.
$$
holds.
\end{lemma}
\begin{proof}
We recall that the constant implicit in the symbol $\ll$ may depend
on the fixed positive integer $k.$ Therefore, we may suppose that
$p>k.$

 Let $\ell$ be an
integer to be chosen later and satisfying the condition
$$
\log p<\ell<p^{\frac{2k+1}{2k+2}}.
$$
Denote by $V$ the set of the first
$\ell$ prime numbers coprime to $\frac{p-1}{d}.$ Since for any
positive integer $m$ there are $O(\log m)$ (even $O(\log m/\log\log
m)$) different prime divisors of $m,$ then from Chebyshev's theorem
it is immediate that for any $v\in V$ we have
$$
v\ll(|\mathcal{V}|+\log p)\log p\ll |\mathcal{V}|\log p,
$$
where $|\mathcal{V}|=\ell$ denotes the cardinality of $\mathcal{V}.$
For a given divisor $d\mid p-1$ denote by $\mathcal{U}_d$ the set of
all elements of the ring $\mathbb{Z}_{(p-1)/d}$ relatively prime to
$(p-1)/d,$ that is $\mathcal{U}_d=\mathbb{Z}_{(p-1)/d}^*.$ For a
given integer $y$ with $(y, (p-1)/d)=1$ consider the congruence
\begin{equation}\label{congruence}
uv\equiv y\pmod{((p-1)/d)}, \quad u\in \mathcal{U}_d,\quad v\in
\mathcal{V}.
\end{equation}
The number of solutions of this congruence is exactly equal to
$|\mathcal{V}|.$ This follows from the fact that once $v$ is fixed
then $u$ is determined uniquely.

We replace $\lambda$ by $g^{t},$ where $t=(p-1)/T,$ and consider the
sum
$$
\sum_{y\in \mathcal{L}_d}\gamma(dy)\mathbf{e}_{p}(ag^{tdxy}).
$$
Let $\delta(y):=\delta(\mathcal{L}_d; y)$ be the characteristic
function of the set $\mathcal{L}_d$ in the ring
$\mathbb{Z}_{(p-1)/d}.$ Since the number of solutions of the
congruence~(\ref{congruence}) is equal to $|\mathcal{V}|$ for any
fixed $y\in \mathcal{L}_d,$ then
$$
\sum_{y\in
\mathcal{L}_d}\gamma(dy)\mathbf{e}_{p}(ag^{tdxy})=\frac{1}{|\mathcal{V}|}\sum_{u\in
\mathcal{U}_d}\sum_{v\in
\mathcal{V}}\gamma(duv)\delta(uv)\mathbf{e}_{p}(ag^{tdxuv}).
$$
Therefore, setting
\begin{equation}
\label{eqn:def} \mathcal{R}_a(\gamma; d, \mathcal{X},
\mathcal{L}_d)=\sum_{x\in \mathcal{X}}\Big|\sum_{y\in
\mathcal{L}_d}\gamma(dy)\mathbf{e}_{p}(ag^{tdxy})\Big|
\end{equation}
we see that
$$
\mathcal{R}_a(\gamma; d, \mathcal{X},
\mathcal{L}_d)=\frac{1}{|\mathcal{V}|}\sum_{x\in
\mathcal{X}}\Big|\sum_{u\in \mathcal{U}_d}\sum_{v\in
\mathcal{V}}\gamma(duv)\delta(uv)\mathbf{e}_{p}(ag^{tdxuv})\Big|,
$$
whence
$$
\mathcal{R}_a(\gamma; d, \mathcal{X}, \mathcal{L}_d)\le
\frac{1}{|\mathcal{V}|}\sum_{u\in \mathcal{U}_d}\sum_{x\in
\mathcal{X}}\Big|\sum_{v\in
\mathcal{V}}\gamma(duv)\delta(uv)\mathbf{e}_{p}(ag^{tdxuv})\Big|.
$$
Application of H\"{o}lder's inequality to the sum over $x$  yields
$$
\mathcal{R}_a(\gamma; d, \mathcal{X}, \mathcal{L}_d)\le
\frac{|\mathcal{X}|^{1-\frac{1}{2k}}}{|\mathcal{V}|} \sum_{u\in
\mathcal{U}_d}\left(\sum_{x=1}^{p-1} \Big|\sum_{v\in
\mathcal{V}}\gamma(duv)\delta(uv)
\mathbf{e}_{p}(ag^{tdxuv})\Big|^{2k}\right)^{1/2k}.
$$
If $(n, p-1)=d,$ if $x$ runs through $\mathbb{Z}_{p-1}$ and if $z$
runs through the reduced residue system modulo $p,$ then $g^{nx}$
and $z^{d}$ run the same system of residues modulo $p$ (including
the multiplicities). Since $(du, p-1)=d,$ then
\begin{equation}
\label{eqn:R2k} \mathcal{R}_a(\gamma; d, \mathcal{X},
\mathcal{L}_d)\le
\frac{|\mathcal{X}|^{1-\frac{1}{2k}}}{|\mathcal{V}|} \sum_{u\in
\mathcal{U}_d}M_u^{1/2k},
\end{equation}
where
$$
M_u=\sum_{z=1}^{p-1} \Big|\sum_{v\in
\mathcal{V}}\gamma(duv)\delta(uv)
\mathbf{e}_{p}(az^{tdv})\Big|^{2k}.
$$
In order to estimate $M_u,$ we observe that
$$
M_u=\sum_{v_1\in
\mathcal{V}}\ldots\sum_{v_{2k}\in\mathcal{V}}\left(\prod_{j=1}^k\prod_{s=k+1}^{2k}
\gamma(duv_j)\overline{\gamma(duv_s)}\delta(uv_j)\delta(uv_s)\right)S(v_1,\ldots,
v_{2k}),
$$
where
$$
S(v_1,\ldots, v_{2k})=\sum_{z=1}^{p-1}
\mathbf{e}_{p}(az^{tdv_1}+\ldots+az^{tdv_k}-az^{tdv_{k+1}}-\ldots
-az^{tdv_{2k}}).
$$
Note that if the set of values $v_1,\ldots, v_k$ is a permutation of
the set of values $v_{k+1},\ldots, v_{2k},$ then
$$
S(v_1,\ldots, v_{2k})=p-1.
$$
Otherwise, we have
$$
az^{tdv_1}+\ldots+az^{tdv_k}-az^{tdv_{k+1}}-\ldots
-az^{tdv_{2k}}=a(a_1z+a_2z^2+\ldots+a_rz^r),
$$
where $a_i, 1\le i\le r,$ are integers and for some $j$ we have
$(aa_j,p)=1$ (here we use that $p>k$). Hence, in this case we can
apply the classical Weil bound (see Chapter~5 of~\cite{LN}) to
obtain
$$
S(v_1,\ldots, v_{2k})\le (td\max_{1\le i\le 2k}v_i)p^{1/2}\ll
td|\mathcal{V}|p^{1/2}\log p.
$$
Therefore, putting all together and using the condition
$|\gamma(n)|\le 1,$ we deduce the bound
$$
M_u\ll\sum_{v_1\in
\mathcal{V}}\ldots\sum_{v_{k}\in\mathcal{V}}\left(\prod_{j=1}^k
\delta(uv_j)\right)p+ \sum_{v_1\in
\mathcal{V}}\ldots\sum_{v_{2k}\in\mathcal{V}}
\left(\prod_{j=1}^{2k}\delta(uv_j)\right)td|\mathcal{V}| p^{1/2}\log
p,
$$
whence
$$
M_u\ll\left(\sum_{v\in\mathcal{V}}\delta(uv)\right)^kp+
\left(\sum_{v\in\mathcal{V}}\delta(uv)\right)^{2k}td|\mathcal{V}|
p^{1/2}\log p.
$$
Combining this estimate with~\eqref{eqn:R2k}, we derive
\begin{eqnarray*}
\mathcal{R}_a(\gamma; d, \mathcal{X}, \mathcal{L}_d)\ll
\frac{|\mathcal{X}|^{1-\frac{1}{2k}}p^{1/2k}}{|\mathcal{V}|}
\sum_{u\in
\mathcal{U}_d}\left(\sum_{v\in\mathcal{V}}\delta(uv)\right)^{1/2}+\\
\frac{|\mathcal{X}|^{1-\frac{1}{2k}}}{|\mathcal{V}|}(td|\mathcal{V}|
p^{1/2}\log p)^{1/2k} \sum_{u\in
\mathcal{U}_d}\sum_{v\in\mathcal{V}}\delta(uv)
\end{eqnarray*}

Since the number of solutions of congruence~(\ref{congruence}) is
equal to $|\mathcal{V}|$ for any $y\in \mathcal{L}_d,$ then
\begin{equation*}
 \sum_{u\in \mathcal{U}_d}\sum_{v\in
 \mathcal{V}}\delta(uv)=|\mathcal{V}||\mathcal{L}_d|.
\end{equation*}
Besides, applying the Cauchy inequality we obtain that
$$
\sum_{u\in \mathcal{U}_d}\left(\sum_{v\in
\mathcal{V}}\delta(uv)\right)^{1/2}\le |U_d|^{1/2}\left(\sum_{u\in
\mathcal{U}_d}\sum_{v\in \mathcal{V}}\delta(uv)\right)^{1/2}\le
(p/d)^{1/2}|\mathcal{V}|^{1/2}|\mathcal{L}_d|^{1/2}.
$$
Therefore,
\begin{equation}
\label{eqn:R2kfin} \mathcal{R}_a(\gamma; d, \mathcal{X},
\mathcal{L}_d)\ll
\frac{|\mathcal{X}|^{1-\frac{1}{2k}}p^{1/2k}p^{1/2}|\mathcal{L}_d|^{1/2}}{|\mathcal{V}|^{1/2}d^{1/2}} +\\
|\mathcal{X}|^{1-\frac{1}{2k}}(td|\mathcal{V}| p^{1/2}\log
p)^{1/2k}|\mathcal{L}_d|.
\end{equation}

Let us define the number of elements of $\mathcal{V}.$ Since
$|\mathcal{L}_d|\le p/d,$  from the condition of the Lemma we see
that
$$
(\log
p)^{3/2}<\frac{p^{\frac{2k+1}{2k+2}}}{d|\mathcal{L}_d|^{\frac{k}{k+1}}t^{\frac{1}{k+1}}(\log
p)^{\frac{1}{k+1}}}<p^{\frac{2k+1}{2k+2}}.
$$
Hence, we can define the number of elements of $\mathcal{V}$ by
taking
$$
\ell=|\mathcal{V}|=\left[\frac{p^{\frac{2k+1}{2k+2}}}{d|\mathcal{L}_d|^{\frac{k}{k+1}}t^{\frac{1}{k+1}}(\log
p)^{\frac{1}{k+1}}}\right].
$$
Inserting this into~\eqref{eqn:R2kfin}, we obtain
$$
\mathcal{R}_a(\gamma; d, \mathcal{X}, \mathcal{L}_d)\ll
\frac{|\mathcal{X}|^{1-\frac{1}{2k}}p^{\frac{1}{2k}+\frac{3}{4k+4}}|\mathcal{L}_d|^{1-\frac{1}{2k+2}}\log
p} {T^{\frac{1}{2k+2}}}.
$$
Lemma~\ref{Lemma} is proved.
\end{proof}

\section{Proof of Theorem \ref{thm:doublesumcomposit}}

If $T^{\frac{2k+1}{2k+2}}\le
|\mathcal{X}|^{1/2k}|\mathcal{Y}|^{-1+\frac{1}{2k+2}}p^{\frac{3}{2}-\frac{1}{2k}-\frac{3}{4k+4}},$
then one can easily check that
$$
\frac{|\mathcal{X}|^{1-\frac{1}{2k}}|\mathcal{Y}|^{1-\frac{1}{2k+2}}
p^{\frac{1}{2k}+\frac{3}{4k+4}}}{T^{\frac{1}{2k+2}}}\ge|\mathcal{X}|^{1-\frac{k+1}{k(2k+1)}}|\mathcal{Y}|
p^{\frac{k+1}{k(2k+1)}}\ge |\mathcal{X}||\mathcal{Y}|.
$$
Hence, in this case the estimate of
Theorem~\ref{thm:doublesumcomposit} becomes trivial. Therefore, we
may suppose that
\begin{equation}
\label{eqn:T}
T^{\frac{2k+1}{2k+2}}\ge
|\mathcal{X}|^{1/2k}|\mathcal{Y}|^{-1+\frac{1}{2k+2}}p^{\frac{3}{2}-\frac{1}{2k}-\frac{3}{4k+4}}.
\end{equation}
Similarly, we may assume that $ T>p^{1/2}(\log p)^{10k}.$

For a given divisor $d\mid p-1$ we denote by $\mathcal{L}_d$ the set
of integers $y$ such that $dy\in \mathcal{Y}$ and $(dy, p-1)=d.$
Then
\begin{equation}\label{eqn:conclusion}
W_{a}(\gamma; T; \mathcal{X}, \mathcal{Y})\le \sum_{d\mid
p-1}\mathcal{R}_a(\gamma; d, \mathcal{X}, \mathcal{L}_d),
\end{equation}
where $ \mathcal{R}_a(\gamma; d, \mathcal{X}, \mathcal{L}_d)$ is
defined by~(\ref{eqn:def}).

Note that $|\mathcal{L}_d|\le \min\{|\mathcal{Y}|, p/d\}.$ For the
divisors $d\mid p-1$ with the condition $d\ge Tp^{-\frac{1}{2}}(\log
p)^{-4k}$ we use the trivial estimate
$$
\mathcal{R}_a(\gamma; d, \mathcal{X}, \mathcal{L}_d)\le
|\mathcal{X}||\mathcal{L}_d|\le
\frac{|\mathcal{X}|p}{Tp^{-\frac{1}{2}}(\log
p)^{-4k}}=\frac{|\mathcal{X}|p^{\frac{3}{2}}(\log p)^{4k}}{T}.
$$
For $d < Tp^{-\frac{1}{2}}(\log p)^{-4k}$ we apply the bound of
Lemma~\ref{Lemma};
\begin{eqnarray*}
\mathcal{R}_a(\gamma; d, \mathcal{X}, \mathcal{L}_d)\ll
\frac{|\mathcal{X}|^{1-\frac{1}{2k}}|\mathcal{Y}|^{1-\frac{1}{2k+2}}p^{\frac{1}{2k}+\frac{3}{4k+4}}\log
p} {T^{\frac{1}{2k+2}}}.
\end{eqnarray*}
Inserting these bounds into~(\ref{eqn:conclusion}) and noting that
$\tau(p-1) = p^{o(1)},$ we deduce the estimate
$$
W_{a}(\gamma; T; \mathcal{X}, \mathcal{Y}) \ll
\frac{|\mathcal{X}|^{1-\frac{1}{2k}}|\mathcal{Y}|^{1-\frac{1}{2k+2}}p^{\frac{1}{2k}+\frac{3}{4k+4}+o(1)}}
{T^{\frac{1}{2k+2}}}+\frac{|\mathcal{X}|p^{\frac{3}{2}+o(1)}}{T},
$$
whence, due to the inequality~\eqref{eqn:T}, we conclude that
$$
W_{a}(\gamma; T; \mathcal{X}, \mathcal{Y})\ll
\frac{|\mathcal{X}|^{1-\frac{1}{2k}}|\mathcal{Y}|^{1-\frac{1}{2k+2}}p^{\frac{1}{2k}+\frac{3}{4k+4}+o(1)}}
{T^{\frac{1}{2k+2}}}.
$$

Theorem~\ref{thm:doublesumcomposit} is proved.



\end{document}